\documentclass[11pt]{article}
\usepackage{enumerate}
\input{amssym.def}
\input{amssym.tex}

 \let\secc=\S

 \def\g{\gamma}                              
 \def\d{\delta}               \def\l{\lambda}              
                                
        \def\n{\nu}                  
 \def\z{\zeta}                \def\x{\xi}

                                                           \def\o{\omega}

                                \def\S{\Sigma}

          \newcommand{\an}{\wedge}

           \def\all{\forall}
           \def\Equi{\Longleftrightarrow}   \def\ex{\exists}

           \def\models{\vDash}
           
           \def\imp{\longrightarrow}
           \def\nat{{\Bbb N}}
        
          \newtheorem{thom}{Theorem}[section]
    \newtheorem{fact}[thom]{Fact}
          \newtheorem{lemma}[thom]{Lemma}
          \newtheorem{propo}[thom]{Proposition}
             \newtheorem{de}[thom]{Definition}
      \newtheorem{corr}[thom]{Corollary}

      \def\thm{\begin{thom}}           \def\ethm{\end{thom}}
          \def\lem{\begin{lemma}}      \def\elem{\end{lemma}}

          \def\prop{\begin{propo}}         \def\eprop{\end{propo}}
          \def\df{\begin{de}}
          \def\edf{\end{de}}
          \def\cor{\begin{corr}}           \def\ecor{\end{corr}}
      \def\fa{\begin{fact}}        \def\ef{\end{fact}}
\def\cala{{\cal A}}

\def\cant{\mbox{}^\omega 2}
\def\defs{ =_{df}}

\def\ie{{\em i.e.}\,}
\def\eg{{\em e.g.}\,}

\def\cf{{\em cf. }\,}

\def\bar{\ove}

\def\la{\langle}
\def\ra{\rangle}

\newcommand{\qd}[1]{\hspace*{\fill{\bf \mbox{ Q.E.D.#1}}}}

\def\pf{{\bf Proof: }}

\def\bar{\overline}

\def\sm{\small}

\def\da{\mathop{\downarrow}}
\def\up{\mathop{\uparrow}}

\begin{document}

\thispagestyle{empty}

\title{ \large\sc $P^f \neq NP^f$ for almost all $f$}
\def\sm{\small}
\author{ J. D. Hamkins\thanks{The research of J. D. Hamkins is supported in part by grants from the NSF and
the PSC-CUNY Research Foundation.}\\
\sm Mathematics, Georgia State University, 30 Pryor Street, Atlanta, GA 30300 USA.\\
\sm Mathematics, City University of New York, 365 Fifth Avenue, New York, NY 10016 USA.\\
\sm http://jdh.hamkins.org \\
       \\
P. D. Welch\thanks{P.D.Welch would like to thank the Deutsche Forschungsgemeinschaft and the University of
Bonn for the support of this research.}\\
\sm Department of Mathematics, University of Bristol, Bristol BS8 1TW, England.\\
\sm  Mathematisches Institut, Beringstra{\ss}e 6, 53115 Bonn, Germany.\\
\sm http://www.maths.bris.ac.uk/$\sim$mapdw\\
}
\date{\today}

\maketitle

\begin{abstract}
{\small \noindent We discuss the question of Ralf-Dieter Schindler whether for infinite time Turing
machines $P^f = NP^f$ can be true for any function $f$ from the reals into $\o_1$. We show that ``almost
everywhere'' the answer is negative. }\end{abstract}

\section{\large\bf{Introduction}}

After establishing $P\not=NP$ for infinite time Turing machines,
Ralf-Dieter Schindler in \cite{S} introduced the more general question
of whether $P^f=NP^f$ for these machines. The classes are defined as
follows:

\eject \df Let $f:\cant \imp \omega_1$ and $A\subseteq \cant $.
\begin{description}\item (A) We say that $A\in P^f$
if there is an (infinite time) Turing machine computable function
$\varphi_e$ so that\\
(i) $A$ is decidable by $\varphi_e$, that is $x \in A$ if and only
if $\varphi_e(x)\da 1$, and\\
(ii) $\forall z \in \cant\quad \varphi_e(z)\da \,\,$ in at most $f(z)$ many steps.
\item (B) We say that $A\in NP^f$ if there is a Turing machine computable
function $\varphi_e$ so that\\
(i) $x \in A$ if and only if there exists $y \in \cant$ so that
$\varphi_e(x\oplus y)\da 1$, and\\
(ii) $\forall z = (x \oplus y) \in \cant\quad \varphi_e(z)\da \,\,$ in at most $f(x)$ many steps.
\end{description}
\edf

The function $f$ plays the role here of the class of polynomials in the classical $P=NP$ question, bounding
the length of the allowed computations. Because of this, one is primarily interested here in the functions
$f$ which are Turing invariant, in the sense that if $x$ and $y$ are Turing equivalent, then $f(x)=f(y)$.
Indeed, since one might expect that a more complicated input should be allowed more time for computation,
it is natural to restrict attention only to the functions $f$ for which $x\leq_T y$ implies $f(x)\leq
f(y)$. The main results of this article, however, do not happen to rely on these assumptions. Since the
computations allow for infinite input, one might usually want to assume that $f(x)$ is infinite.

If the value of $f(x)$ is some constant $\alpha$ then the classes $P^f$
lie strictly within the Borel hierarchy (\cite{S} Lemma 2.7). If $f(x)
= \o_1^x$ then $P^f$ coincides with hyperarithmetic (and so we are
really still within the realms of Kleene recursion \eg see \cite{HS}).
When $f(x) > \o_1^x$ for all $x$ we then truly enter for the first time
the world of sets that are essentially computed by infinite time Turing
machines (see \cite{HL} for the basic concepts). \cite{S} raises two
questions concerning these classes for such functions dominating $f(x)
= \o_1^x$.

Note that unlike the basic notions of $P$ and $NP$
Schindler defined in \secc 2 of \cite{S}, his class
$NP^f$ is not, in general, just the projection of
$P^f$.

We wish to prove that for almost all functions $f$ the classes $P^f$ and $NP^f$ are different. Given the
extra information that the verifying witness $y$ can encode, this is, or should be, unsurprising. The first
point to note is that if the values $f$ takes are sufficiently large, they will exceed the times needed by
a machine to establish membership of any decidable set.

We recall a definition from \cite{HL}: \df $\lambda^x \defs
\sup\{\alpha \mid \ex e\, \varphi_e(x)\da y \an y \in WO \an rk(y) =
\alpha \}$. \edf

Equivalently (and the reader may take this as a definition):

\fa \label{lambda}(\cite{W6} Theorem 1.1) $\lambda^x$ is the supremum
of halting times of any Turing computable function on input $x$. \ef

 Implicit in this latter result---when taken with the definition
of decidable sets of reals \cite{HL}---(see the discussion in
\cite{W8}) is the following characterisation of such sets.

\fa $A \in \cant$ is decidable if and only if there  are $\S_1$ formulae in the language of set theory
$\varphi_0(v_0),\varphi_1(v_0)$ so that $$ x \in A\quad \Equi\quad L_{\l^x}[x]\models
\varphi_0[x]\quad\Equi\quad L_{\l^x}[x]\models\neg\varphi_1[x] $$ \ef

Clearly then:

\lem If $\l^x\leq f(x)$ for all $x$, then $P^f$ is the class of
decidable sets. \elem

We recall a definition from \cite{W6}: \df \label{sigma} $\Sigma^x
\defs  \sup\{\,\alpha\, \mid \mbox{ a
code $x \in WO$ for $\alpha$ occurs on a tape }$\\
of some computation $\varphi_e(x)$ at some time $\}$ \edf

(Thus $\Sigma^x$ is the supremum of ordinals with so-called {\em
accidentally writable} (relative to $x$) codes, as defined in \cite{HL}
p. 580.)
 It is shown in \cite{W} that (i) $\S^x$ is the least
$\sigma$ so that
 $L_\sigma[x]$ has a proper $\Sigma_2$-elementary substructure: and in
\cite{W6} that
 (ii)  $L_{\Sigma^x}[x]$ is inadmissible.

Our first theorem is that any function $f$ that dominates the function
$g(x) = \S^x$ separates the classes $NP^f$ and $P^f$. This answers
Question 2 of \cite{S} ``almost everywhere":

\thm \label{every} Let $f$ satisfy $f(x) \geq \Sigma^x $ for all $x$.
Then $NP^f \supsetneqq P^f$. \ethm

The next theorem answers Question 1 of \cite{S}. Let $f_0$ be defined
as $f_0(x) = \o_1^x +\o $.

\thm \label{adm}
 $NP^{f_0}
\supsetneqq P^{f_0}$. \ethm

In \cite{S} the classes $P^{f_0}$ and $NP^{f_0}$ are denoted $P^{++}$ and $NP^{++}$, respectively, so we
have shown $P^{++}\neq NP^{++}$, the question left open in \cite{S}. As \cite{S} points out, of course
$P^{f_0}$ contains all $\S^1_1$ and $\Pi^1_1$ sets. It is not hard to see that if $A \in \mathrm{Diff}( <
\o_1^{ck}, \S^1_1)$ the Hausdorff difference hierarchy for levels below the first non-recursive ordinal,
then $A \in P^{f_0}$.

Variants of these methods  will also show that many other functions whose range falls between these two
ordinals will also separate these two classes.  We have not attempted an exhaustive classification.

Suppose $\phi_e: \cant \imp \nat$ is total. Then the length of the
computation $\phi_e(x)$, $\n_x$ say, defines for us  a ``clock''.
Namely let $f(x) =
\n_x$; suppose for convenience $\n_x$ is always of limit length and at least $\o^\o$.

\thm \label{nclock} With $f$ as above:
 $NP^{f}
\supsetneqq P^{f}$. \ethm

\section{Preliminaries}

We shall let $\omega_1^x$ stand for the first ordinal not recursive in
$x$. Then $L_{\omega^x_1}[x]$ is an {\em admissible } set. We refer the
reader to \cite{Bar} for an account of admissible sets and their basic
properties. We shall use the following notation for the machine
configurations. Let the cells of the tape be enumerated $\la C_i|i <
\o\ra$ with the cell $C_i$ having value  $C_i(\xi)$ at time $\xi$. We
assume that the first $n$ blocks on the tape are enumerated by $\la
C_i|i < 3n\ra$ with $C_0,C_1,C_2$ being the leftmost output,scratch,
and input cells respectively. A {\em snapshot} of the tape at time
$\gamma$ is then a function $s \in \cant$ coding these cell values,
with $s(i) = C_i(\gamma)$ (possibly also allowing it to encode
somewhere internal states, the location of the head and the instruction
of the program about to be performed). A halting computation is then
entirely given by the wellordered snapshot sequence of computations of
the length of the computation. The machine is considered to be
specified by a finite program, just as for ordinary Turing machines,
although the head is allowed to read, and write to, triples of cells at
any one stage. Thus a typical instruction might be of the form $(q_i,
j, j', X, q_k)$ where $j,j' \in \mbox{}^32$, interpreted to mean that
in state i viewing cells $C_{3l},C_{3l},C_{3l+2}$ with values
$j(0),j(1),j(2)$ the machine moves to state $q_k$, changes the cell
values to those of $j'$ and moves one unit in the direction $X \in
\{L,R\}$ (for Left and Right). The machine has however a special limit
state $q_L$ (and at limit times it is in this state viewing
$C_{0},C_{1},C_{2}$). The machine may thus halt at a limit time if it
contains a quadruple of the form $(q_L,  j, j', q_H)$. Note however
that executing this last step of computation means that it changes the
entries of $C_{0},C_{1},C_{2}$ to those given by $j'$. It is this
feature that allows the classes $P^{f_0}$ to be closed under
complementation: an ``accepting'' entry of $C_0$ as 1 can be switched
to a ``rejecting'' zero at the last moment. To be completely clear, if
the machine executes the halt instruction at the limit stage $\n$, we
reiterate that the final value is $C_i(\nu+1)=j'(i)$, which might
differ from the limit value $C_i(\nu) = j(i)$.

We shall use the following fact.

\fa \label{canon} \cite{Ga} There is an index $e_0 \in \o$ so that,
uniformly, for any $x \in\cant$ $\{e_0\}^x$ is an illfounded linear
ordering of $\o$, recursive in $x$, with wellfounded part of order type
$\o^x_1$. \ef

We use this index to give us a ``canonical $\o_1^x + \omega$-clock'':
an algorithm that halts in exactly $\o_1^x + \omega$ steps. (The
following argument is the ``uniform in $x$'' version of that of Theorem
3.2 of \cite{HL}.) The algorithm does the following: it first computes
the field of the relation $ <_R = \{e_0\}^x$ and then proceeds, by
picking the least element of the field, $n$ say,
 in $\omega$ many steps to find  the $<_R$-least element below $n$
of the ordering. It then, in another $\omega$ many
steps, proceeds to strike out all mention of this
element from the field of $ \{e_0\}^x$; it then
picks an element $n'$ of the field that is left and
searches for the next $<_R$-least element below
this $n'$; at limit stages below $\o_1^x$ the
procedure continues smoothly as the wellfounded
part of $<_R $ has order type $\o_1^x$. However in
the interval $( \o_1^x, \o_1^x + \o)$ it searches
in vain for a least element. It chooses some $\bar
n$ in the field that is left. We may assume that
each time it descends in $<_R$ it flashes a signal
in $C_2$, by alternating in the next 3 stages the
value of $C_2$ to be ``0,1,0''. After $\omega$ many
stages, by the limsup rule of the machines, the
value $C_2(\o_1^x + \o) = 1$, and moreover this is
the first time this happens at a limit stage. We
assume then the program has been written so as to
immediately halt if this occurs.

An alternative $\omega_1^x+\omega$ clock is obtained by the algorithm
that on input $x$ simulates all computations on input $x$, looking for
a stage at which none of the programs halt. Since $\omega_1^x$ is the
least such stage, and it takes $\omega$ many steps to recognize that
this situation has occurred, the algorithm can halt exactly at
$\omega_1^x+\omega$.

\df\label{eventual} $\zeta^x \defs  \sup\{\,\alpha \mid
\ex e \ex y \ex \g \all \d >\g \mbox{ $y \in WO$  lies on the output tape }$\\
$\varphi_e(x)$ at time $\d \an rk(y) = \alpha\,\}.$
\edf

We shall appeal to the following fact:

\fa \label{lzs} (``$(\l,\z,\S)$ Theorem'')\begin{description}\item
(i)(\cite{W6} \cf 2.3.) For any computation of the form $\phi_e(x)$,
the snapshot at time $\z$, $s_\z$, is exactly that at time $\S^x$,
$s_\S$; they are both {\em settled} snapshots, \ie they are destined to
recur on a closed and unbounded class of ordinals;
\item (ii) (\cite{W} 2.1,2.3) if  $ \x $ is the least $\xi$ satisfying
$L_\xi[x]\prec_{\Sigma_2} L_{\Sigma^x}[x]$, then
(a) $\x = \z^x$ and (b) $(\l^x ,\z^x, \S^x)$ is the
lexicographic least increasing triple $(\l,\z,\S)$
satisfying $L_\l[x]\prec_{\Sigma_1}
L_\z[x]\prec_{\Sigma_2}
L_{\Sigma}[x]$.\end{description} \ef

By \cite{W6} (Claim (ii) of 3.4) there are computations $\phi_{q_0}$ so
that for any input $x$, $(\z^x,\S^x)$ is the lexicographically (on
$On\times On$) least pair of ordinals with repeating snapshots
$(s_{\z^x},s_{\S^x})$: running a universal machine provides such. (In
fact, any computation on input $x$ which does not repeat before
$\lambda^x$ is such an example. In such a case it is the snapshot
$s_{\z^x}$ that provides a parameter witnessing the inadmissability of
$L_{\S^x}[x]$ - \cf \cite{W6}, 3.4) In general then, this pair of
snapshots witnesses that the computation is either halted or in an
infinite loop.

\section{Separating the classes}



{\bf Proof of Theorem \ref{every}}. By our observation in the
introduction, under these assumptions $P^f$ is the class of decidable
sets of reals. Let $ H = \{ \la p,x\ra \in \o \times \cant \mid
\phi_p(x)\da \,\,\, \}$. Then $H$ is the complete set coding the
halting problem for sets of reals. $H$ is undecidable, but the above
arguments, together with the $(\l,\z,\S)$ Theorem will show that $H \in
NP^f$.

It will suffice to verify whether $\phi_e(x)\da\,\,\,$ by the following
method.

We consider informally a Turing machine computable algorithm $P_n$ that
effects the following:

$P_n$ on input $e\smallfrown x\oplus y$: \\
(i) First checks whether $y$ codes an $\omega$-model containing (an
isomorphic copy of) $x$ and: $$ \cala = \langle \o,y\rangle \models ``
KP \an V = L[x]\an  \S^x \mbox{ exists} \an  \phi_e(x)\da\,\,\,."$$

By way of explanation: we intend that $\cala $ thinks there is a least
initial segment of its $L[x]$-hierarchy with a proper $\S_2$ elementary
substructure - this is the import of ``$\S^x$ exists." This is an
arithmetic condition on $x \oplus y$ and thus can be checked by $P_n$
in $< \o^\o$ many steps. If this fails for $y$ then $P_n(e \smallfrown
x\oplus y)\da 0$ thus halting with a zero in the first cell $C_0$ of
the output tape.

(ii) Otherwise a preliminary ``1" is written to
$C_0$ and then $P_n$ proceeds to check if
$WFP(\cala)$ contains the true $ \S^x$. However we
first dispose of the part of the model containing
all sets of $(L[x]$-rank$)^\cala \geq
(\S^x)^\cala$. We simply eliminate all reference to
these in $y$, thus in effect rewriting $y$ as some
new real $\bar y$. However this is again a simple
operation, and can be done in, say, $\omega^2$
steps (note that $y$ has some integer $n$ which
denotes the $(\S^x)^\cala$ so it is a trivial
matter to do this). The process then proceeds to
check for the wellfoundedness of the ordinals of
this new $\o$-model $\bar \cala = \la \o,\bar y\ra$
determined by the initial segment $(\S^x)^\cala$ in
the usual way by erasing integers from the field of
its ordinals.

If the model $\bar \cala$ is wellfounded then this process takes the true $\S^x$ many steps, (note that
$(\S^x)^\cala \cong \S^x$ as there is $n$ denoting $(\S^x)^\cala$ and the property of an ordinal being
$\S^x$ is absolute), and furthermore the model $\bar\cala$ is correct about $\phi_e(x)\da$. Using the trick
of keeping track of when the least (in some standard ordering of $\o \times \o$) pair is erased we may
realise also that the field of the ordering of $On^{\bar \cala}$ has become empty (\cf the proof of Theorem
3.1 in \cite{HL}). If so it can halt exactly at the $\S^x$'th step with the required $1$ in $C_0$.

If the model $\bar \cala$ is  illfounded (and hence
$(\Sigma^x)^\cala$ is in the illfounded part of the
original $On^\cala$), then in fact  $\S^x \nsubseteq
 WFP(\bar\cala)=WFP(\cala)$:
this is because (a) we cannot have $\S^x \in WFP(\bar\cala) $ (as otherwise $\cala$ would recognize it as
$\S^x$); (b) but neither can $WFP(\bar\cala) = \S^x$ (as  $L_{\S^x}[x]$ is inadmissible, and this would
contradict the Truncation Lemma (\cite{Bar}).

Hence any instance of illfoundedness in
$(On)^{\bar\cala}$ will be detected before the true
$\S^x$ many steps have been taken. This leaves time
to change the contents of $C_0$ to a zero, and halt
- here before $\S^x$ many steps have been taken.

In each case then $P_n(e\smallfrown x\oplus y)$ halts in no more than
$\S^x$ many steps with the correct output.

\qd{(Theorem \ref{every})}

The algorithm above can be made more time efficient, so that confirming instances of the decision problem
are settled more quickly. This modified algorithm can be made to actually follow the naive idea that to
determine whether $\phi_e(x)\da$, one should simply simulate the computation $\phi_e(x)$ to see if it
halts, and somehow end simulations that have gone on too long. The point is that the model-checking method
of the previous argument, where one checks whether $\bar\cala$ is well-founded, is essentially a
nondeterministic clock for $\Sigma^x$, in the sense that it halts at time $\Sigma^x$ for certain witnesses
$y$, and before $\Sigma^x$ for all other witnesses. Our modified algorithm, therefore, is simply to run
such a clock alongside the computation of $\phi_e(x)$, and accept the input if the computation halts before
the clock runs out. Since in the worst case the clock runs to time $\Sigma^x$, this algorithm
nondeterministically decides whether $\phi_e(x)\da$ in time $\Sigma^x$. But the point is that affirmative
instances are decided much earlier, in time before $\l^x$, because this is when the halting computations
actually halt.

The theorem can be improved by ignoring the bold-face context of the situation:

\thm Suppose that $f(p)\geq\Sigma$ for every finite $p$ and $f(x)\geq\omega$ for all other reals $x$.  Then
$NP^f\supsetneqq P^f$.\label{lightface}\ethm

\pf The idea is that the (weak) halting problem $h=\{\, p\, \mid
\varphi_p(0)\da\,\, \}$ will be in $NP^f$ but not in $P^f$. It clearly
is not in $P^f$, since it is not decidable. But one can see it is in
$NP^f$ by the following algorithm:  on input $x\oplus y$, first check
whether $x$ is a finite $p$ or not. If not, then halt and reject the
input. Otherwise, carry out the algorithm of Theorem \ref{every} on the
input $p\smallfrown 0\oplus y$. With suitable choice of $y$, this will
decide whether $(p,0)\in H$, which is equivalent to $p\in h$, in at
most $\Sigma^0=\Sigma$ many steps, as desired.\qd{}\\

{\bf Proof of Theorem \ref{nclock}} Let $H_{f} = \{ \la p,x\ra \in \o
\times \cant \mid \phi_p(x)\da \mbox{ in } \leq \n_x \mbox{ steps } \}$

(1) $H_{f} \notin P^{f}$.

\pf Let $r$ be the partial function defined as follows:
$$r(y) = \left\{
\begin{array}{cl}{0}
 & \mbox{ if } y = \la p,x \ra \an y \notin H_{f} \\
 \up & \mbox{ otherwise}
\end{array}
\right.
$$
If ``$y \in H_{f}$'' were decidable by an algorithm that always halted
in at most $ \n_y $ steps then $r$ would also be computable by an
algorithm $P_m$, that if it converged, would do so in   $ \n_y $ steps.
(We could obtain a program for $r$ by simply changing the behaviour of
that of the former algorithm by switching at the very last limit step
where it halted on a 1, into some non-halting loop.) Let $\phi_m: \cant
\imp \o$ be this latter function. Let $c_0$ be the constant zero
function. However then
$$
\begin{array}{lcr}
\la m, c_0\ra \in H_{f} &\Equi &
 \phi_m( \la m, c_0\ra) \da
\mbox{ in } \leq \n_x \mbox{ steps }\\  & \Equi & r( \la m, c_0\ra ) =
0  \Equi  \la m, c_0\ra \notin H_{f}
\end{array}
$$
a contradiction. \qd{(1)}

(2)   $H_{f} \in NP^{f}$.

\pf We use the ideas from the proof of Theorem \ref{every}. We devise
an algorithm $P_n$ to verify (2). $P_n$ on input $p\smallfrown x\oplus
y$:

 (i) First checks whether $y$ codes an
 $\omega$-model containing (an isomorphic copy of)
 $x$ and:
 $$ \cala = \langle \o,y\rangle \models `` KP \an V =
 L[x] \an  \phi_p(x)\da \mbox{ in }
 \leq \phi_e(x) = \n_x \mbox{ steps}."$$

This is again an arithmetic condition on $x \oplus y$ and thus can be
checked by $P_n$ in $< \o^\o$ many steps. If this fails for $y$ then
$P_n(p \smallfrown x\oplus y)\da 0$.

(ii) Otherwise a preliminary ``1" is written to
$C_0$ and then $P_n$ proceeds to check if
$WFP(\cala)$ contains the true $ \n_x$. As before
we dispose of the part of the model containing all
sets of $(L[x]$-rank$)^\cala \geq (\n^x)^\cala$.
The process then proceeds to check for the
wellfoundedness of this new initial segment model
$\bar \cala$ up to $On^{\bar \cala} \cong
(\n_x)^\cala$. We use that \cite{HL} (Theorem 8.8)
shows  $\n_x$ is not an $x$-admissible ordinal.

If the model is wellfounded then this process takes the true $\n_x$
many steps,  thus it can halt exactly at the $\n_x$'th step with the
required $1$ in $C_0$.

%
Arguing as before using the cited inadmissibility
of $L_{\n_x}[x]$, the wellfounded part of $\cala$
cannot have rank exactly the true $\n_x$; hence we
are justified in testing only the initial segment
of the ordinals of $\bar \cala$ determined by
$(\n_x)^\cala$. Then any instance of illfoundedness
will be encountered strictly before $\n_x$ many
steps have been taken.

In each case then $P_n(p\smallfrown x\oplus y)$ halts in no more than
$\n_x$ many steps with the correct output.

\qd{(Theorem \ref{every})}

As a final comment some of the above discussion may lead one to
considering the class of sets $A$ such that $x\in A$ and $x\notin A$
can each be verified quickly, that is, such that there are two
programs, such that $x\in A$ if and only if there is a witness $y$ such
that $x\oplus y$ is accepted by the first program, and $x\notin A$ if
and only if there is a witness $y$ such that $x\oplus y$ is accepted by
the second program, and both programs halt on any input $x\oplus y$ in
time before $f(x)$ if they halt at all. That is, both $x\in A$ and
$x\notin A$ can be verified quickly, with the correct choice of
verifying witnesses, but there is no insistence that the programs
compute quickly (or even halt at all) when given irrelevant witnesses
verifying nothing. B. L\"{o}we has pointed out that such a class of
sets corresponds to a notion of  $ NPTIME^f \cap co$-$NPTIME^f$, but we
have made no investigation of such concepts.


\bibliographystyle{amsplain}
\bibliography{settheory2}

\end{document}